\documentclass[12pt]{amsart}
\usepackage{amsmath}
\usepackage{amstext}
\usepackage{amsbsy}
\usepackage{amssymb}
\usepackage{amsfonts}
\usepackage{latexsym}
\usepackage{verbatim}

\newfont{\script}{eusb10}
\numberwithin{equation}{section}
\newcommand{\C}{\mathbb C}
\newcommand{\R}{\mathbb R}

\newcommand{\Z}{\mathbb Z}

\newcommand{\s}{\mbox{\boldmath $\sigma$}}
\newcommand{\CP}{\C \mathbb{P}^1}
\newcommand{\MT}{\mbox{\small{$\widetilde{M}$}}}
\newcommand{\M}{\mathcal{M}}

\newcommand{\Id}{\mathrm{Id}}

\newcommand{\res}{\xi_{\mbox{\tiny{$-1$}}}}

\newcommand{\SL}{\mathrm{SL}(2,\C)}
\newcommand{\Sl}{\mathfrak{sl}(2,\C)}

\newcommand{\SU}{\mathrm{SU}(2)}
\newcommand{\su}{\mathfrak{su}(2)}

\newcommand{\gl}{\mathfrak{gl}(2,\C)}

\newcommand{\LSU}{\Lambda  \SU_\sigma}

\newcommand{\LGE}{\Lambda_r \SU_\sigma }
\newcommand{\LGC}{\Lambda_r \SL_\sigma}
\newcommand{\LGI}{\Lambda_r^+ \SL_\sigma}

\newcommand{\LAC}{\Lambda_r \Sl_\sigma}

\newcommand{\Lsu}{\Lambda \su_\sigma}

\newcommand{\pott}{\Lambda \Omega(\MT)}

\newcommand{\Lsl}{\Lambda_{-1}^\infty \Sl_\sigma}

\newcommand{\be}{\begin{equation}} 
\newcommand{\ee}{\end{equation}}

\newtheorem{theorem}{Theorem}[section]

\newtheorem{lemma}[theorem]{Lemma}

\title{On the associated family of Delaunay surfaces}
\author{M. Kilian}

\address{Martin Kilian, 
Department of Mathematical Sciences,
University of Bath, Bath, BA2 7AY, United Kingdom}
\email{masmk@maths.bath.ac.uk}

\thanks{ Mathematics Subject 
Classification 53A10.}

\begin{document}

\begin{abstract}
We use the method of Dorfmeister, Pedit and Wu 
\cite{DorPW} to obtain the associate 
family of Delaunay surfaces and derive a 
formula for the neck size of the 
surface in terms of the entries of the holomorphic 
potential.
\end{abstract}
\maketitle

\section{Introduction} 
The Gauss map of a constant mean curvature (CMC) 
surface is a harmonic map \cite{RuhV} 
and Uhlenbeck \cite{Uhl} has shown that such 
harmonic maps can be obtained as 
projections of horizontal holomorphic maps from the 
universal cover of the surface into a certain loop 
group. There is a Weierstra\ss\ type representation 
for this procedure, commonly referred to as the 
DPW method, due to 
Dorfmeister, Pedit \& Wu \cite{DorPW}. 
The construction involves solving a linear 
differential system with values in 
a loop group and then Iwasawa decomposing the 
solution to obtain the extended unitary frame of 
the Gauss map.

The DPW method is analogous to the 
Weierstra\ss\ representation for minimal 
surfaces insomuch as input data in the form 
of meromorphic functions on a Riemann surface 
generate a conformal CMC immmersion of the 
universal cover. If the underlying Riemann 
surface is not simply connected, 
it is necessary to investigate various monodromies 
and solve period problems in order that one member 
of the resulting associate family is invariant under 
deck transformations. Recently, progress has 
been made in constructing CMC immersions of 
$n$--punctured spheres with Delaunay ends in 
the DPW approach \cite{KilMS}, \cite{KilSS} 
and \cite{Sch}. 

Since Delaunay surfaces are ubiquitous in the study 
of CMC surfaces with embedded ends \cite{KorKS}, 
we feel the necessity to give a self contained 
account of Delaunay surfaces in the DPW framework.  
Here we thus desribe the construction of 
the associated family of Delaunay 
surfaces using holomorphic data on the twice 
punctured Riemann sphere. 

This paper grew out of a part of the author's 
PhD thesis \cite{Kil} and it is a pleasure to 
herewith express my gratitude to my advisor 
Franz Pedit as well as to Nicholas Schmitt 
for numerous useful discussions.

\section{Loop groups}
We adopt the notation 
$\mathrm{diag}[u,\, v] = 
\left( \begin{smallmatrix} 
u & 0 \\ 0 & v \end{smallmatrix} \right)$, 
$\mathrm{off}[u,\, v] = 
\left( \begin{smallmatrix} 
0 & u \\ v & 0 \end{smallmatrix} \right)$ and 
begin by collecting some well known 
results on loop groups.

For real $r \in (0,1]$ denote the analytic maps of 
$C_r = \{ \lambda \in \C : | \lambda | = r \}$ 
with values in $\SL$ by 
$\Lambda_r \SL = \mathcal{O}(C_r,\SL)$. 
We have an involution on maps $C_r \to \gl$
given by
\begin{equation*} 
  \s : g(\lambda) \mapsto 
  \sigma \, g(-\lambda)\, \sigma^{-1} \mbox{ with }\,
  \sigma = \mathrm{diag}[1,-1]
\end{equation*}
and denote the \textit{twisted} $r$--Loop group 
of $\SL$ by
\begin{equation*}
  \LGC = \left\{ g \in \Lambda_r \SL: 
  \, \s g = g \right\}.
\end{equation*}
The Lie algebras of these groups, denoted by
$\LAC$, consist of analytic maps
$g:C_r\rightarrow \Sl$, which satisfy 
$\s g = g$. We will use the 
following subgroups of $\LGC$:
Let $K = \left\{ \mathrm{diag}[a,1/a] : 
a > 0 \right\} \subset \mathrm{SL}(2,\R)$ and 
$I_r = \{ \lambda \in \C : |\lambda|<r \}$ 
and denote
\begin{equation*}
  \LGI = \left\{ g \in \LGC \cap 
  \mathcal{O}(I_r,\SL) : 
  g(0) \in K \right\}.
\end{equation*}
Let $A_r = \{ \lambda \in \C : 
r<|\lambda|<1/r \}$ and by an abuse of 
notation set
\begin{equation*}
  \LGE = \left\{ g \in \LGC \cap 
  \mathcal{O}(A_r,\SL): 
  \left. g \right|_{S^1} \in \SU \right\}.
\end{equation*}
For $r=1$ we will omit the subscript.
Corresponding to these subgroups, we analogously 
define Lie subalgebras of $\LAC$. Of fundamental 
importance is a certain Loop group factorisation, 
the Iwasawa decomposition. 
We cite McIntosh \cite{McI}, where it is proven 
that Multiplication 
\begin{equation*}
  \LGE \times \LGI \to \LGC
\end{equation*}
is a diffeomorphism onto. The associated splitting 
$g = FB$ of $g \in \LGC$ with $F \in \LGE$ and 
$B \in \LGI$ will be called Iwasawa decomposition. 
The condition $B(0) \in K$ 
ensures that the factorization is unique and 
\begin{equation}\label{eq:Iwaunique}
\LGE \cap \LGI = \left\{ \Id \right\}.
\end{equation}
%

%%%%%%%%%%%%%%%%%%%%%%%%%%%%%%%%%%%%%%%%%

\section{DPW method}

Let $\Omega(M)$ denote the holomorphic $1$--forms on 
a Riemann surface $M$ and define 
\begin{equation*}
  \Lsl = \left\{ \xi \in 
  \mathcal{O}(\C^*,\Sl_\sigma): 
  \lambda \xi \in \mathcal{O}(\C,\Sl) \right\}.
\end{equation*}
CMC surfaces come in $S^1$ 
families, the \textit {associated family}. 
The DPW method \cite{DorPW} constructs 
conformal CMC immersions of the universal cover 
$\MT$ in the following three steps. 
Let 
\begin{equation*}
  \xi \in \pott = \Omega(\MT) \otimes \Lsl,
\end{equation*}
$\tilde{z}_0 \in \MT$ and $\Phi_0 \in \LGC$: 
The first step consists in solving the 
initial value problem 
\begin{equation}\label{eq:IVP}
  d\Phi = \Phi \xi,\, \Phi(\tilde{z}_0)= \Phi_0
\end{equation}
to obtain a unique map $\Phi : \MT \to \LGC$. 
Secondly, Iwasawa decompose $\Phi = F \, B$ 
pointwise on $\MT$ to obtain a unique map 
$F : \MT \to \LGE$. \\
Let $\partial_\lambda = 
\tfrac{\partial}{\partial \lambda}$. 
The final step is performed by plugging 
$F$ into the Sym-Bobenko formula 
\begin{equation} \label{eq:Sym}
  f_\lambda = -\tfrac{1}{H} \left(i \lambda 
  (\partial_\lambda F) \,F^{-1} + \tfrac{i}{2} 
  F \,\sigma\, F^{-1} \right) 
\end{equation}
to obtain conformal immersions 
$f_\lambda :\MT \to \Lsu$ with 
constant mean curvature $H \neq 0$, 
that is for each $\lambda_0 \in S^1$ 
we have a conformal CMC immersion 
$f_{\lambda_0} :\MT \to \su \cong \R^3$. 
We call $\xi \in \pott$ the 
\textit{holomorphic potential}. 
We thus have a map $(\xi,\Phi_0,\tilde{z}_0) 
\mapsto f_\lambda$ and the triple 
$(\xi, \, \Phi_0, \, \tilde{z}_0)$ is called the 
DPW data of a CMC immersion. 
The metric and Hopf differential are determined by 
a triple $(\xi, \, \Phi_0, \, \tilde{z}_0)$ as 
follows. Locally, on $U \subset M$, 
if we write 
$\xi = \lambda^{-1}\mathrm{off}[a_1,a_2] + \ldots$ 
and $B(\lambda =0) =\mathrm{diag}[r,r^{-1}]$ 
for holomorphic $1$--forms 
$a_i\in\Omega(U)$ and 
$r:U \to \R^+$, then it can be shown 
that $f_1$ has metric $4r^4|a_1|^2$ and Hopf 
differential $Q=-\tfrac{1}{2}a_1a_2$. It follows 
that $f_1$ has branch points at the zeroes of 
$a_1$ and umbilics at the zeroes of $a_2$. 
Our basis of $\su$ will be 
\begin{equation*}
  e_{1}=\begin{pmatrix} i&0\\0&-i \end{pmatrix},\,  
  e_{2}=\begin{pmatrix} 0&1\\-1&0 \end{pmatrix},\, 
  e_{3}=\begin{pmatrix} 0&i\\i&0 \end{pmatrix}.
\end{equation*}
As a consequence of using twisted loop groups 
in the construction, there is some redundency 
in the associated family in that every 
surface appears twice up to rigid motion. 
\begin{lemma}
Under the map $\lambda \mapsto -\lambda$ the 
corresponding surface is rotated about $e_1$ 
by $\pi$ - radians, that is 
$f_\lambda \mapsto e_1 f_\lambda e_1^{-1}$.
\begin{proof}
Twistedness is $\xi (-\lambda ) = 
\s \xi (\lambda ) = e_1 \xi(\lambda) e_1^{-1}$. 
If $\Phi$ solves \eqref{eq:IVP} with triple 
$(\xi,\Phi_0 ,\tilde{z}_0 )$, then $\Phi e_1$ 
solves \eqref{eq:IVP} with triple 
$(e_1 \xi e_1^{-1},\Phi_0 e_1 ,\tilde{z}_0 )$.
\end{proof}
\end{lemma}
It is a well known, that to every CMC surface 
there exist two parallel surfaces, one of constant 
Gaussian curvature, the other of constant mean 
curvature, obtained by moving every point on the 
surface in the normal direction by $\tfrac{1}{2H}$ 
and $\tfrac{1}{H}$ respectively. The parallel 
associated family of CMC immersions is given by 
\begin{equation}
  f_{\lambda}^{\mbox\tiny{par}} = -\tfrac{1}{H} 
  (i\lambda (\partial_\lambda F) F^{-1} - 
  \tfrac{i}{2} F\,\sigma \,F^{-1}).
\end{equation}
\begin{lemma} \label{th:parallel}
Let $f_\lambda$ be the associated family of 
CMC immersions obtained from 
$(\xi , \Phi_0 ,\tilde{z}_0 )$. Then the triple 
$(e_3 \xi e_3^{-1} ,e_3\Phi_0e_3^{-1},\tilde{z}_0 )$ 
generates the parallel associated family 
$f_\lambda^{\mbox\tiny{par}}$.
\begin{proof}
If $\Phi$ solves \eqref{eq:IVP} for the triple 
$(\xi,\Phi_0 ,\tilde{z}_0 )$, then 
$\Psi = e_3\Phi e_3^{-1}$ solves \eqref{eq:IVP} 
for the triple 
$(e_3 \xi e_3^{-1} ,e_3\Phi_0e_3^{-1},\tilde{z}_0 )$ 
with resulting associated family $g_\lambda$. 
If $\Phi=FB$ is the Iwasawa decomposition of $\Phi$, 
then $\Psi=e_3 FB e_3 ^{-1}$. Since $b_0$ is 
diagonal, we can write $b_0 e_3 ^{-1}  =  
e_3 ^{-1} b_0 ^{-1}$ and thus $B e_3 ^{-1} =  
e_3 ^{-1} \widetilde{B}$ for some 
$\widetilde{B} \in \LGI$. Hence 
$\Psi=e_3F \widetilde{B}e_3^{-1}$ and the 
Sym-Bobenko formula gives
\begin{align*}
g_\lambda  &= -\tfrac{1}{H}(i\lambda e_3 
(\partial_\lambda F) F^{-1} e_3 ^{-1} + 
\tfrac{1}{2} e_3 F e_3 ^{-1} e_{1} e_3 
F^{-1} e_3 ^{-1} ) \\
&= -\tfrac{1}{H} e_3(i\lambda (\partial_\lambda F) 
F^{-1} -\tfrac{1}{2} F e_{1} F^{-1})e_3^{-1} 
= e_3 f_\lambda^{\mbox\tiny{par}} e_3^{-1} 
\end{align*}
Therefore the surface generated by the potential 
$e_3\xi e_3^{-1}$ is the parallel surface rotated 
about $e_3$ by $\pi$ - radians. 
\end{proof}
\end{lemma}
The left action by $\LGE$ on the initial condition, 
$\Phi_0 \mapsto C \Phi_0$ for $C \in \LGE$ results 
in an Euclidean motion of the surface. 
In fact, we have the following   
\begin{lemma}\label{th:rigidmotion}
Let $\Phi$ solve \eqref{eq:IVP} with triple 
$(\xi, \Phi_0,\tilde{z}_0 )$ and let 
$f_\lambda$ denote the resulting associated 
family. Let $C \in \LGE$ and for $\alpha \in \C$ 
with $|\alpha|=1$ define 
$A = \rm{diag}[\alpha,\bar{\alpha}]$. 
Then $\Psi=C \Phi A$ solves \eqref{eq:IVP} 
for the triple 
$( A^{-1} \xi A , C \Phi_0 A,\tilde{z}_0 )$ 
and the resulting associated family $g_\lambda$ 
is is given by  
\begin{equation}
  g_\lambda = -\tfrac{1}{H} i \lambda 
  (\partial_\lambda C) C^{-1} + C f_\lambda C^{-1}.
\end{equation}
\begin{proof}
If $\Phi=FB$ is the Iwasawa decomposition of 
$\Phi$, then $A^{-1} B A \in \LGI$ and the 
Iwasawa decomposition of $\Psi$ is 
$\Psi = CFBA = CFAA^{-1}BA$. Hence $\Psi$ 
has the unitary frame $CFA$, which inserted 
into the Sym--Bobenko formula gives
\begin{align*}
  g_\lambda 	&= -\tfrac{1}{H} 
  (i\lambda (\partial_\lambda CFA) 
  (CFA)^{-1} + \tfrac{1}{2} (CFA)e_{1}(CFA)^{-1}) \\
  &= -\tfrac{1}{H} (i\lambda ( 
  (\partial_\lambda C) C^{-1} + 
  C (\partial_\lambda F )
  F^{-1}C^{-1} +  \tfrac{1}{2} 
  CFe_{1}F^{-1}C^{-1} )) \\
  &= -\tfrac{1}{H} i \lambda (\partial_\lambda C ) 
  C^{-1} + C f_\lambda C^{-1}
\end{align*}
an proves the result. 
\end{proof}
\end{lemma}
%

%%%%%%%%%%%%%%%%%%%%%%%%%%%%%%%%%%%%%%%%%%%%
%%%%%%%%%%%%%%%%%%%%%%%%%%%%%%%%%%%%%%%%%%%%

\section{The period problem}

\noindent Let $M$ be a connected Riemann surface 
with universal cover $\MT$ and $\Delta$ 
the group of deck transformations. Let 
$\xi \in \pott$ with $\gamma^*\xi = \xi$ for all 
$\gamma \in \Delta$. Let $\Phi: \MT \to \LGC$ 
be a solution of the differential equation 
$d \Phi = \Phi \xi$. Writing 
$\gamma^* \Phi = \Phi \circ \gamma$ for 
$\gamma \in \Delta$, we define 
$\chi(\gamma) \in \LGC$ by 
$\chi(\gamma) = (\gamma^* \Phi) \, \Phi^{-1}$. 
The matrix $\chi(\gamma)$ is called the 
\textit{monodromy} matrix of $\Phi$ with respect 
to $\gamma$. 

If $\widehat{\Phi}: \MT \to \LGC$ 
is another solution 
of $d \Phi = \Phi \xi$ and 
$\widehat{\chi}(\gamma)=
(\gamma^*\widehat{\Phi})\,\widehat{\Phi}^{-1}$, 
then there exists a constant $C \in \LGC$ such that 
$\widehat{\Phi} = C \Phi$. Hence
$\widehat{\chi}(\gamma) = C \chi(\gamma) C^{-1}$ and 
different soultions give rise to 
mutually conjugate monodromy matrices. 

A choice of base point 
$\tilde{z}_0 \in \MT$ and initial condition 
$\Phi_0 \in \LGC$ gives 
the \textit{monodromy representation} 
$\chi : \Delta \to \LGC$
of $\xi$. 
Henceforth, when we speak of the monodromy 
representation, or simply monodromy, 
we tacitly assume that it is induced 
by an underlying triple 
$(\xi, \Phi_0, \tilde{z}_0 )$ with 
invariant holomorphic potential 
$\gamma^*\xi = \xi$ for all 
$\gamma \in \Delta$. It is shown in 
\cite{DorH:cyl} that CMC immersions of open 
Riemann surfaces $M$ can always be generated by such 
invariant holomorphic potentials. 

If $\Phi = F B$ is the pointwise Iwasawa 
decomposition of $\Phi : \MT \to \LGC$, 
then we shall need to study the monodromy of $F$ 
as a means to controlling the periodicity of the 
resulting immersion \eqref{eq:Sym}. 
A priori, we are not assured that 
the quantity $\M(\gamma) = ( \gamma^*{F} ) 
F^{-1}$ is $z$--independent for all 
$\gamma \in \Delta$. One way to circumvent 
this issue is to ensure that $\chi$ is 
$\LGE$--valued. 
\begin{lemma}\label{th:Fholonomy}
Let $\chi(\tau)$ be the monodromy of $\xi$ 
for $\tau \in \Delta$. If $\chi(\tau) \in \LGE$ 
then $\chi(\tau) = \M(\tau)$.
\begin{proof}
From the Iwasawa decomposition $\Phi = FB$ we obtain 
$\chi(\tau) F B = \tau^*F \tau^*B$. 
If $\chi(\tau) \in \LGE$, then 
\eqref{eq:Iwaunique} yields
$(\tau^*F)^{-1} \chi(\tau) F  = 
(\tau^*B) B^{-1} = \mathrm{Id}$. 
Therefore $\chi(\tau) = \M(\tau)$.
\end{proof}
\end{lemma}

The condition $\chi(\tau) \in \LGE$ can be ensured 
if the initial condition $\Phi_0$ in \eqref{eq:IVP} 
is unitary and the potential $\xi$ is skew hermitian 
along a curve $\tau \in \pi_1(M)$ that passes 
through the base point.
\begin{lemma}\label{th:skewhermitian}
  Let $\xi \in \pott$ and $\tau^* \xi = \xi$ for 
  $\tau \in \pi_1(M) \cong \Delta$. 
  Pick a point $z_0 \in \tau$ and 
  $\Phi_0 \in \LGE$ and let 
  $\Phi:\MT \to \in \LGC$ 
  be the solution of \eqref{eq:IVP} with triple 
  $(\xi, \Phi_0, z_0)$. If  $\xi$ is $\Lsu$--valued 
  along $\tau$, then $\chi (\tau) \in \LGE$.
\begin{proof}
The potential $\xi$ being $\Lsu$--valued along 
$\tau$ means that $\bar{\xi}^t = -\xi$ along 
$\tau$. Hence for $\xi = \Phi^{-1}d\Phi$ we have 
$d\Phi \bar{\Phi}^t + \Phi d\bar{\Phi}^t = 0$. 
Integrating this last equation gives 
$\Phi \bar{\Phi}^t = C$ for a $z$--independent 
$C \in \LGC$. The initial condition satisfies 
$\Phi_0 \overline{\Phi_0}^t = \mathrm{Id}$ forcing 
$C =\mathrm{Id}$. Hence $\Phi$ is $\LGE$--valued 
along $\tau$ ensuring that $\chi (\tau) \in \LGE$. 
\end{proof}
\end{lemma}    
The next result characterizes 
the period problem in the DPW framework. 
\begin{theorem}\cite{DorH:per}\label{th:closeperiod1}
Let $\xi \in \pott$ have $\LGE$--valued monodromy 
$\chi$ and let $f_\lambda$ be the associated family 
generated by the triple 
$(\xi, \rm{Id},\tilde{z}_0)$. 
Then there exists a $\lambda_0 \in S^1$ such that 
$\tau^* f_{\lambda_0} = f_{\lambda_0}$ for 
$\tau \in \Delta$ if and only if 
$\chi(\tau,\,\lambda)$ satisfies both
\begin{align}
  \chi (\tau,\,\lambda_0) &= 
  	\pm \rm{Id},\label{eq:closing1}\\
  \left.\partial_{\lambda} \chi (\tau) 
  \right|_{\lambda_0} &= 0 
\label{eq:closing2}.
\end{align}
\end{theorem}
%

%%%%%%%%%%%%%%%%%%%%%%%%%%%%%%%%%%%%%%%%%%%%%%%%%%%%%
%%%%%%%%%%%%%%%%%%%%%%%%%%%%%%%%%%%%%%%%%%%%%%%%%%%%%

\section{Delaunay surfaces}

As an application of the above, 
we shall construct the holomorphic 
data in the DPW representation that generates the 
family of Delaunay surfaces. The conformal structure 
of a Delaunay surface is that of a twice punctured 
Riemann sphere, so we let $M=\C^*$. If we 
identify $\C^* \cong \C / 2\pi i \Z$, then 
$\exp:\C \to \C^*$ is the universal 
covering map and we shall write 
$\tilde{z} = \log(z)$. The group of deck 
transformations $\Delta \cong \Z$ is generated by 
$\tau: \tilde{z} \mapsto \tilde{z} + 2 \pi i$.  
In the associated family, we shall solve the period 
problem \eqref{eq:closing1} and \eqref{eq:closing2} 
for $\lambda_0=1$, and choose the base point 
$\tilde{z}_0=0$. Consider 
\begin{equation} \label{eq:Delaunay}
  \xi = \res \frac{dz}{z} \mbox{ with }
  \res = \begin{pmatrix} 
   c & a \lambda^{-1} + \bar{b}\lambda \\
   \lambda^{-1} + \bar{a} \lambda & -c \end{pmatrix}
\end{equation}
with $ a,\,b \in \C$ and $c \in \R$. 
The two eigenvalues of $\res$ are 
$\pm \sqrt{ - \det \res}$. Since the residue 
is invariant under $z \mapsto \alpha z$ for 
$|\alpha|=1$, the induced surface has intrinsic 
rotational symmetry. The solution to 
\eqref{eq:IVP} with triple 
$(\xi,\rm{Id},0 )$ is obtained by 
exponentiation and with $z=re^{it}$ can be factored 
\begin{equation}
  \Phi (z,\lambda)= \exp ( \log (z) \, \res ) =
  \exp ( it \, \res ) \exp ( \ln |r| \, \res ).
\end{equation}
\begin{lemma}\label{th:framerot}
$F_1 (t) := \exp (i\,t\,\res ) \in \LSU$ is a 
one parameter family of rotations with 
common axis $i\res$.
\begin{proof}
Since $\res$ is hermitian 
$it\res \in \mathfrak{su}_2$ for all $t \in \R$.
Hence $F_1 (t) \in \LSU$ and $[ i\res, F_1 (t) ]=0$. The claim now follows by uniqueness of the axis.
\end{proof}
\end{lemma}
\noindent The differential equation for a profile 
curve $r(t)$ of a CMC $H$ surface of revolution is 
given by 
\begin{equation}
  \frac{\ddot{r} + 1/r}{(1 + \dot{r}^2)^{3/2}} = 2H
\end{equation}
which is equivalent to the first order equation 
$\dot{s}=0$ for 
$s(t) = \tfrac{r}{\sqrt{1 + \dot{r}^2}} - r^2 H$. 
At a neck or a bulge we have $\dot{r} = 0$, so there 
is a constant $\kappa$ such that $r(1-rH)=\kappa$. 
In 1841 Delaunay \cite{Del} proved the following 
theorem. An elementary proof is given by 
Smyth \cite{Smy}.
\begin{theorem} \cite{Del} A surface of revolution 
is a constant mean curvature surface if and only 
if the profile curve is a roulette of a conic.
\end{theorem}
The resulting one parameter family of 
CMC cylinders, complete immersed surfaces 
of revolution, are called Delaunay surfaces. 
A Delaunay surface is uniquely determined 
by the bulge and neck radii up to rigid motions. 
The surfaces obtained from ellipses are called 
unduloidal, the surfaces obtained from parabolas 
are called nodoidal \cite{Eel}. 
In particular, Delaunay 
surfaces have no umbilic point. A quadratic 
meromorphic differential on $\CP$ has degree $-4$. 
For an unbranched cylinder without umbilics, 
the four poles of the Hopf differential have 
to be at the ends. By translational symmetry, 
the Hopf differential of a Delaunay surface has 
to have a pole of order $2$ at each end. Hence, 
up to coordinate transformations, the Hopf 
differential for a Delaunay surface is of the 
form $z^{-2}dz^2$. 
\begin{theorem}\label{th:Delaunay}
Let $\xi = \res\tfrac{dz}{z}$ whith $\res$ as in 
\eqref{eq:Delaunay}. If
$|a + \bar{b}|^2 + c^2 = 1/4$ and $ab \in \R$, 
then the triple $(\xi,\Id,0)$ generates the 
associated family of Delaunay surfaces which 
at $\lambda = 1$ has bulge/neck radius
\begin{equation} \label{eq:radius}
  \tfrac{1}{2H}(\,1\pm\sqrt{1 - 16ab}\,).
\end{equation}
\begin{proof}
Diagonalizing $\res$, the solution to \eqref{eq:IVP} 
with triple $(\xi,\Id,0)$ is given by
\begin{equation*}
  \Phi(z,\lambda) = T(\lambda) \mathrm{diag}
  [z^{-\mu(\lambda)},z^{\mu(\lambda)}] 
  T(\lambda)^{-1} \mbox{ where }T(\lambda) = \left( 
  \begin{smallmatrix} 
  \tfrac{(c - \mu)\lambda}{ \bar{a}\lambda^2 + b} & 
  \tfrac{(c + \mu)\lambda}{ \bar{a}\lambda^2 + b} \\
  1 & 1 \end{smallmatrix} \right)
\end{equation*}
and $\mu(\lambda) = \sqrt{ a \overline{a} + 
b \overline{b} + c^2 + ab\lambda^{-2} + 
\overline{ab}\lambda^2 }$. Since 
$\xi \in \Lambda^{\sigma} \mathfrak{su}_2$ 
along $|z|=1$, the holonomy is $\LGE$--valued and 
given by 
\begin{equation*}
  \chi(\lambda) = T(\lambda) \,\mathrm{diag}[
  e^{-2\pi i\mu(\lambda)},\,e^{2\pi i\mu(\lambda)}] 
  \, T(\lambda)^{-1}.
\end{equation*}
The first closing condition $\chi(1) = \pm \Id$ is 
equivalent to $\mu (1) \in \tfrac{1}{2}\Z$. For the 
constants $a,b,c$ this holds if and only if 
$|a + \bar{b}|^2 + c^2 \in \tfrac{1}{4} \Z$. 
The second closing condition 
$\partial_\lambda H |_{\lambda=1} = 0$ holds 
if $\partial_\lambda \mu |_{\lambda=1} = 0$ 
which is equivalent to $ab \in \R$ if the 
first closing condition is assumed. 
Let $\arg$ denote some fixed 
branch of argument, then since $ab \in \R$, 
we have $\arg(\bar{b})=\arg(a)$. 
Let $\gamma=e^{i\arg(a)}$ and define 
$A=\mathrm{diag}[ \sqrt{\gamma}, 
\sqrt{\bar{\gamma}}]$. 
By Lemma \ref{th:rigidmotion}, the surface generated 
by the triple $(A \xi A^{-1}, \Id, 1)$ 
is a rotation of the surface generated by the triple 
$(\xi,\Id,1)$, but since the coefficients of 
$\lambda^{-1},\ \lambda$ in the 
potential $A \xi A^{-1}$ are real, 
we can assume without loss of generality 
that $a,b \in \R$ and the two closing 
conditions together with the requirement 
that the surface be simply wrapped, 
reduce to the equation of an elliptic 
cylinder $( a + b )^2 + c^2 = \tfrac{1}{4}$. 
This leaves two real paramters of freedom, 
the necksize of the surface and the radius 
of the image circle of $S^1$: At the 
basepoint $z_0=1$ we have 
$\Phi(1,\lambda)=\Id$ and the image 
under the immersion is 
$f(1,\lambda) = -\tfrac{1}{2H} e_1$. 
By Lemma \ref{th:framerot}, the frame along $|z|=1$ 
is given by the one parameter family of rotations 
$\exp(it \res)$ with common axis 
$i\res$. Hence the resulting surface is a 
surface of revolution and thus a Delaunay surface. 
The unit cirle $|z|=1$ is mapped to the circle 
with center 
\begin{equation*}
  \tfrac{1}{2}(f(1,\lambda)- f(-1,\lambda)) =  
\frac{8a^2 + 8ab -1}{2H}\,e_1 -\frac{4ac}{H}\,e_3,
\end{equation*}
radius $2|a\,H^{-1}|$ and axial direction 
$i\res$. The angle $\theta$ between the axis 
and the tangent vector of the profile curve 
$r(t)$ at some $r_0$ is given by 
\begin{equation*}
  \cos \theta = \frac{1}{\sqrt{1 + \dot{r}^2}} = 
  2(a+b).
\end{equation*}
The quadratic equation at a neck or bulge is 
\begin{equation}
  r(1-rH)=2(a+b) 2|a/H| - 4H|a/H|^2
\end{equation}
whose solutions are given by \eqref{eq:radius}. 
\end{proof}
\end{theorem}
By Lemma \ref{th:parallel}, the parallel CMC 
surface of a Delaunay surface is obtained by 
interchanging neck and bulge.  
We conclude our discussion of Delaunay surfaces 
by indicating how the choice of the constants 
determine the resulting surface. 
The case $ab>0$ results in 
unduloids, while $ab<0$ yields nodoids. 
For $a=b=1/4$, the resulting surface is the round 
cylinder and the limiting case $a=1/2,\,b=0$ 
produces a round sphere with two points removed. 
%
%%%%%%%%%%%%%%%%%%

\bibliographystyle{amsplain}
\bibliography{ref}

\end{document}